\documentstyle{amsppt}
\NoBlackBoxes
\topmatter
\title
A note on hamiltonian Lie group actions and Massey products
\endtitle
\author Zofia St\c epie\'n and Aleksy Tralle
\endauthor
\address
Szczecin Technical University, al. Piast\'ow 38, 70-323 Szczecin, Poland
\vskip6pt
University of Warmia and Mazury, ul. \D Zo\l\/nierska 14A, 10-561 Olsztyn, Poland
\endaddress
\email
stepien\@arcadia.tuniv.szczecin.pl,
tralle\@matman.uwm.edu.pl
\endemail
\subjclass
53C15, 55P62
\endsubjclass
\abstract
In this note we show that the property of having only vanishing triple Massey products in the equvariant cohomology is inherited by the set of fixed points of hamiltonian circle actions on closed symplectic manifolds. This result can be considered in a more general context of characterizing homotopic properties of Lie group actions. In particular it can be viewed as a partial answer to the Allday-Puppe question about finding conditions ensuring the "formality" of $G$-actions. 
\endabstract
\endtopmatter
 
\document
\head 1. Introduction
\endhead
In this article we prove the following theorem.
\proclaim{Theorem 1.1} Let $(M,\omega)$ be a closed symplectic manifold endowed with a hamiltonian circle action $S^1\times M\to M$. Denote by $F$ the connected component of the fixed point set of this action. If there exists a non-vanishing triple Massey product in $H^*(F)$, then the same is valid for the equivariant cohomology $H^*_{S^1}(M)$.
\endproclaim
Of course the above theorem can be rephrased as follows.
\proclaim{Corollary 1.1} If $(M,\omega)$ is a closed symplectic manifold endowed with a hamiltonian circle action $S^1\times M\to M$ such that all the triple Massey products in the equivariant cohomology $H^*_{S^1}(M)$ vanish, then the connected component of the fixed point set $F=M^{S^1}$ also has only vanishing Massey products in $H^*(F)$.
\endproclaim
Thus, the vanishing of the triple Massey products is inherited by the connected component of the fixed point set. This result can be considered in a more general context of the cohomology theory of Lie transformation groups \cite{AP}. Note that Allday and Puppe asked a question about possible characterizations of "formal" Lie group actions. Recall that a topological space $X$ is called formal if it has minimal model $\Cal M_X$  quasi-isomorphic to the cohomology algebra $H^*(\Cal M_X)\cong H^*(X)$ (we don't give any details here, referring to  \cite{AP,TO}). It is natural to ask if one can characterize (in terms of the invariants of the given $G$-action)  group actions with the formal fixed point set $X^G$. Since the vanishing of Massey products is an obstruction to formality \cite{TO} (see also references therein),  Theorem 1.1 is a partial answer to the above question. It appears, that the weaker property of vanishing of the triple Massey products is indeed inherited by the fixed points of the {\it hamiltonian} circle actions.
\vskip6pt
On the other hand, this result can be viewed as a characterization of the hamiltonian circle action itself, and in this respect develops the theory along the lines of Kirwan's result \cite{K}. It was proved in \cite{K} that for a closed  $(M,\omega)$ endowed with a hamiltonian  action of a compact Lie group $G$ the Leray-Serre spectral sequence associated with the Borel fibration
$$M\to EG\times_G M\to BG$$
collapses at $E_2$. Note that characterizing hamiltonian and non-hamiltonian actions is a very important question in symplectic topology (cf. \cite{Au, G,JK,LO,McD, McDS}). For example, it is well known, that in dimension 4, a circle action is hamiltonian if and only if it has fixed points \cite{Au, McD}. This fact does not hold in higher dimensions, but the counterexample  of McDuff is quite subtle. It is still not known if there exist non-hamiltonian circle actions with {\it discrete} set of fixed points \cite{G}.  
\subhead Acknowledgment
\endsubhead
We are grateful to Chris Allday who attracted out attention to the problem of "formality" of Lie group actions and shared with us his insights in the theory of group actions.
\vskip6pt
The second author acknowledges the support of the Polish Research Committeee (KBN).
\head 2. Preliminaries and notation
\endhead

\subhead  Massey products and the Cartan model
\endsubhead
Here we recall briefly the notion of Massey products in the form suitable for our considerations. To get a more detailed exposition of this topic in this form, one can consult \cite{RT}. Note however, that we use the term "non-vanishing" or "nonzero" Massey product instead of "essential" Massey product in \cite{RT}.
Let there be given a commutative differential graded algebra $(A,d)$. The cohomology algebra of $(A,d)$ is denoted by $H^*(A)$. If $a\in A$ is a cocycle, we write $[a]$ for the corresponding cohomology class. For a homogeneous element $a\in A$ of degree $p$ we write $\bar a$ denoting $\bar a=(-1)^pa$. Assume that we are given a triple of the cohomology classes $[a],[b],[c]$ such that $[a][b]=[b][c]=0$. Consider $x\in A$ and $y\in A$ such that $dx=\bar ab$ and $dy=\bar bc$. One can check that the element $\bar ay+\bar xc$ is a cocycle and therefore determines the cohomology class $[\bar ay+\bar xc]\in H^*(A)$. Note that this class depends on the choice of $x$ and $y$. By definition the {\it set} of all cohomology classes $[\bar ay+\bar xc]$ is denoted by $\langle [a],[b],[c]\rangle$ and is called the {\it triple Massey product}. 
\definition{Definition} We say that the set $\langle [a],[b],[c]\rangle$ is {\it defined} if $[a][b]=[b][c]=0$ and {\it does not vanish} if the set of all cohomology classes $[\bar ay+\bar xc]$ does not contain zero. In the opposite case we say that $\langle [a],[b],[c]\rangle$ vanishes.
\enddefinition
One can easily describe the indeterminacy in the definition of the triple Massey product and formulate it in the following way (see \cite{RT, Prop. 1.5}). Let $\langle [a],[b],[c]\rangle$ be a defined triple Massey product in $H^*(A)$. Denote by $([a],[c])$ the ideal in $H^*(A)$ generated by elements $[a]$ and $[c]$. The product $\langle [a],[b],[c]\rangle$ does not vanish if and only of there exists the cohomology class $x\in\langle [a],[b],[c]\rangle$ such that $x\not\in ([a],[c])$.
\vskip6pt
We need also the following formulas (see \cite{RT, Prop. 1.4}):
$$
\aligned
\xi\langle a_1,a_2,a_3\rangle &\subset \langle \xi a_1,a_2,a_3\rangle\\
\xi\langle a_1,a_2,a_3\rangle &\subset \langle a_1,\xi a_2,a_3\rangle\\
\xi\langle a_1,a_2,a_3\rangle &\subset \langle a_1,a_2,\xi a_3\rangle
\endaligned
\eqno (2.1)
$$
which are valid for any $a_1,a_2,a_3\in H^*(A)$ and for any $\xi$ represented by a central element in $A$. 
\vskip6pt
Let $f: (A,d)\to (A',d')$ be a morphism of differential graded algebras. Then
$$f^*\langle [a],[b],[c]\rangle\subset \langle f^*[a],f^*[b],f^*[c]\rangle \eqno (2.2)$$
(see Prop. 1.3 in \cite{RT}).
\vskip6pt
Consider now the case of a $G$-manifold, i.e. a manifold endowed with a smooth action of a Lie group.
We use the {\it equivariant} cohomology of the $G$-manifold, i.e. the cohomology of the total space of the Borel fibration: 
$$M\to EG\times_GM\to BG$$
associated with the universal principal $G$-bundle $G\to EG\to BG$ over the classifying space $BG$ of the Lie group $G$. Thus, $H_G^*(M)=H^*(EG\times_GM)$ (see \cite{Au, GS}.
\vskip6pt
In the proof of  Theorem 1.1 we will calculate Massey products with respect to the {\it Cartan model}. Recall that for any $G$-manifold $M$ one can associate the following differential graded algebra. Consider $\Omega^*_G(M)=(\Omega^*(M)\otimes S(\frak g^*))^G$, where $\Omega^*(M)$ is the de Rham algebra of $M$, and $S(\frak g^*)$ is a symmetric algebra over the dual to the Lie algebra $\frak g$ of $G$. Then $G$ acts on $\frak g^*$ by the coadjoint representation and hence there is a natural $G$-action on the tensor product $\Omega^*(M)\otimes S(\frak g^*)$. We consider the subalgebra $\Omega^*_G(M)$ of the fixed points of the given action. The details of this construction can be found in \cite{BV, GS, JK, McDS}. We use the fact that there is a natural differential $D: \Omega^*_G(M)\to \Omega^*_G(M)$ and that 
$$H^*(\Omega_G^*(M),D)\cong H^*G(M).$$
\subhead Hamiltonian $G$-actions
\endsubhead
Finally, recall the definition of the hamiltonian action of a Lie group on a  symplectic manifold \cite{Au, GS, K, McDS}. Let $(M,\omega)$ be a symplectic manifold and let $G$ be a  Lie group acting on $M$ by symplectomorphisms. A smooth map
$$\mu: M\to\frak g^*$$
is called a {\it moment map} of the given $G$-action, if
\roster
\item $\mu$ is $G$-equivariant with repect to the given $G$-action on $M$ and with respect to the coadjoint $G$-action on $\frak g^*$;
\item for any $m\in M$, $X\in\frak g$ and $v\in T_mM$ the following equality holds
$$T_m\mu(v)(X)=\omega(v,\tilde X_m)$$
(where $\tilde X$ denotes the fundamental vector field on $M$ generated by $X$).
\endroster
\definition{Definition 2.2} We call a symplectic $G$-action on $(M,\omega)$ {\it hamiltonian} if it has a moment map.
\enddefinition 
 
The moment map need not be unique. It is not very difficult to check that any symplectic action of a compact Lie group on a closed symplectic manifold $M$ such that $H_1(M)=0$, is hamiltonian. On the other hand, in case of manifolds with large fundamental groups, the notion of hamiltonian action is related to interesting and non-trivial topological questions (cf. \cite{LO,McD}).
\head 3. Proof of  Theorem 1.1
\endhead   
Consider the triple $(M,\omega, G)$. Recall that if the $G$-action is symplectic, the fixed point set $M^G$ is  a symplectic submanifold. Let $F=M^G_0$ be the connected component of $M^G$. We have a symplectic embedding
$$i: F\to M.$$
We can consider $F$ as a symplectic $G$-manifold with a trivial action of $G$. In particular, $H^*_G(F)=H^*(F)\otimes S(\frak g^*)^G$. In case $G=S^1$ we have
$$H^*_{S^1}(F)=H^*(F)\otimes \Bbb R[h] \eqno (3.1)$$
where $\Bbb R[h]$ denotes the free polynomial algebra with one generator $h$ of degree 2. Passing to the Borel fibrations we can write the following commutative diagram.
$$
\CD
F @>i>> M\\
@VVV @VVV\\
EG\times_GF @>{i_G}>> EG\times_GM\\
@VVV @VVV\\
BG @>{=}>> BG
\endCD
$$
On the cohomology level we will get the maps $i^*_G: H^*_G(M)\to H_G^*(F)$ and $i^*: H^*(M)\to H^*(F)$. Let $\nu$ denote the normal bundle of the symplectic embedding $i$. Since $G$ acts on this bundle fiberwise, one can define the {\it equivariant} normal bundle $\nu_G=EG\times_G\nu$. In this way we obtain the vector bundle
$$EG\times_G\nu\to E_G\times_GF.\eqno (3.2)$$
In particular, there is a well defined notion of the {\it equivariant Euler class} of $\nu$. The cohomology class $\chi\in H^*_G(F)$ is called the equivariant Euler class of $\nu$, if $\chi$ is the Euler class of the bundle $(3.2)$.
\vskip6pt
The following facts can be found in \cite{McDS, pp. 192-193} or in \cite{JK} and can be summarized in the following proposition.
\proclaim{Proposition 3.1} Let $G=S^1$ act on a closed symplectic manifold $(M,\omega)$ and let $F=M^G_0$. Then:
\roster
\item "(i)" the normal bundle $\nu$ splits into the sum of complex line bundles
$$\nu=\oplus_j^mL_j,\,\,m=\text{codim}_MF$$
invariant with respect to $S^1$-action. The circle group acts on each $L_j$ with weight $k_j$,
\item "(ii)" the equivariant Euler class has the form
$$\chi=\prod_{j=1}^m(c_1(L_j)+k_jh).\eqno (3.3)$$
\endroster
\endproclaim
The following properties of the map $i_G^*$ characterize the {\it hamiltonian} $G$-action.
\proclaim{Proposition 3.2} \cite{Au,p.139} Assume that $G$ is a torus and acts on a closed symplectic manifold $(M,\omega)$ in a hamiltonian way. Then
\roster
\item "(i)" the map $i_G^*: H^*_G(M)\to H^*_G(F)$ is injective;
\item "(ii)" there exists a linear map (the Gysin homomorphism) $(i_G)_*: H^*_G(F)\to H^*_G(M)$ with the property:
$$i^*_G(i_G)_*(x)=\chi x\,\,\text{for any}\,\,x\in H^*_G(F).$$
\item "(iii)" $\chi$ is not a zero divisor in $H^*_G(F)$.
\endroster
\endproclaim
\proclaim{Corollary 3.1} Under the assumptions of Proposition 3.2 the Gysin map $(i_G)_*$ is injective.
\endproclaim
\demo{Proof} Suppose $x\in\operatorname{Ker}\,(i_G)_*$. If $x\not=0$ Proposition 3.2 (ii) will give $i^*_G(i_G)_*(x)=\chi x=0$. Since $i^*_G$ is injective by Proposition 3.2, the latter equality will give $\chi x=0$, a contradiction with Proposition 2 (iii).
\hfill$\square$
\enddemo
The proof of the Theorem will now follow from two lemmas below. Note that in both lemmas we keep the same notation, and we assume that $G$ is a {\it torus} acting on $(M,\omega)$ in a hamiltonian way.
\proclaim{Lemma 3.1} Let there be given a triple Massey product $\langle u,v,w\rangle\subset H^*_G(F)$. Then there is a defined Massey product
$$\langle u\chi,v\chi,w\chi\rangle. \eqno (3.4)$$
If Massey product (3.4) does not contain zero, then the expression
$$\langle (i_G)_*u,(i_G)^*v,(i_G)^*w\rangle\subset H^*_G(M) \eqno (3.5)$$
defines a non-vanishing Massey product in $H^*_G(M)$.
\endproclaim
\demo{Proof} To avoid clumsy notation let us denote $(i_G)_*$ as $i_*$ and $i^*_G$ as $i^*$.
To show that (3.4) is defined one can write the following equalities
$$i^*(i_*ui_*v)=i^*i_*u\cdot i^*i_*v=\chi u\cdot \chi w=0.$$
Obviuosly, the same could be written for $i_*vi_*w$. Note that the Gysin map is not multiplicative, but $i^*$ is, and this allows us to complete the proof. The following formulae show that (3.4) does not contain zero.
$$i^*\langle i_*u,i_*v,i_*w\rangle\subset \langle i^*i_*u,i^*i_*v,i^*i_*w\rangle =\langle u\chi,v\chi,w\chi\rangle.$$
Here we used (2.1) and Proposition 3.2.
\hfill$\square$
\enddemo
\proclaim{Lemma 3.2} Assume $G=S^1$ and that $\langle u,v,w\rangle$ is a non-vanishing triple Massey product in $H^*(F)\subset H^*_{S^1}(F)=H^*(F)\otimes \Bbb R[h]$. Then 

$$\langle\chi u,\chi v,\chi w\rangle\not=0.$$
\endproclaim
\demo{Proof} According to Section 2, we will prove that there exists an element $z$ in the set $\langle \chi u,\chi v, \chi w\rangle$ such that $z\not\in (\chi u,\chi w)$. Use the equivariant cohomology $H^*_{S^1}(F)$ calculated with respect to the Cartan model. Take a non-trivial Massey triple product $\langle u,v,w\rangle$ considered as a non-trivial Massey product in the equivariant cohomology (one can easily check by straightforward calculation that $\langle u,v,w\rangle$ cannot become zero in the tensor product $H^*(F)\otimes \Bbb R[h]$ writing the corresponding cocycles in the Cartan model). From Lemma 3.1, $\langle \chi u,\chi v,\chi w\rangle$ is defined.  Since $\langle u,v,w\rangle\not=0$, there exists $x\in \langle u,v,w\rangle$ such that $x\not\in (u,w)$. Note that $\chi^3 x\in\langle \chi u,\chi v,\chi u\rangle$ (by (2.2)). Assume that
$$\langle \chi u,\chi v,\chi w\rangle=0.$$
It means that  any $z\in\langle \chi u,\chi v,\chi w\rangle$ belongs to the ideal $(\chi u,\chi w)\subset H^*_{S^1}(F)$. In particular, $\chi^3x\in(\chi u,\chi w)$. Hence
$$\chi^3x=\chi ua+\chi wb,\,a,b\in H^*_{S^1}(F).$$
Therefore
$$\chi(\chi^2x-ua-wb)=0.$$
Recalling that $\chi$ is not a zero divisor (Proposition 3.2 (iii)) one can write
$$\chi^2x=ua+wb.$$
Taking into consideration that $u,w\in H^*(F)\subset H^*_{S^1}(F)$ and the expression for the Euler class (3.2) one obtains
$$\prod_{j=1}^m(c_1(L_j)+k_jh)^2x=u(a_0+a_1h+...+a_{2m}h^{2m})+w(b_0+b_1h+...+b_{2m}h^{2m}).$$
Since $h$ is a free generator, comparing the coefficients yields $x\in (u,w)$, a contradiction.
Finally, $z=\chi^3x$ is the required element.
\hfill$\square$
\enddemo

Now we can complete the proof of  Theorem 1.1. If $\langle u,v,w\rangle$ is a non-trivial Massey triple product in $H^*(F)$, Lemma 3.2 implies that $\langle\chi u,\chi v,\chi w\rangle$ is a non-vanishing triple Massey product in $H^*_{S^1}(F)$. From Lemma 3.1 we get a nontrivial triple Massey product expressed by formula (3.4).
\hfill$\square$
\remark{Remark} In \cite{McD} an example of non-hamiltonian circle action on some closed 6-dimensional manifold was given. Theorem 1.1 may yield the other method of constructing circle actions of this kind. Indeed, if one constructed a symplectic $G$-manifold $M$ with vanishing triple Massey products but with $F=M^G_0$ having a non-vanishing one, this would give an example of non-hamiltonian action having fixed points. However, this does not look very easy.
\endremark

\enddocument